\documentclass[11pt]{amsart}

\usepackage{tikz}
\usepackage{amsmath}
\usepackage{amsthm}
\usepackage{amsfonts}
\usepackage{amssymb}
\usepackage[all]{xy}
\usepackage{alltt}
\usepackage{enumitem}
\usepackage{stmaryrd}
\usepackage{comment}
\usepackage{xspace}
\usepackage{comment}
\usepackage{libertine}
\usepackage[libertine]{newtxmath}
\usepackage{graphicx}
\theoremstyle{plain}

\theoremstyle{definition}

\newcommand{\rar}{\longrightarrow}

\usepackage[margin=1.25in]{geometry}

\usepackage[
        colorlinks=true,urlcolor=blue,linkcolor=blue,citecolor=blue
]{hyperref}

\begin{document}  
 
\title{Rabbit Hunting using Set Theory and Probability}
\date{\today}
 
\author{Sunil Chebolu}
\address{Department of Mathematics \\
Illinois State University \\
Normal, IL 61790, USA}
\email{schebol@ilstu.edu}

\author{Deepayan Sarkar} 
\address {Theoretical Statistics and Mathematics Unit\\
Indian Statistical Institute, Delhi Centre \\
New Delhi 110016, India}
\email{deepayan@isid.ac.in}


 

 
\maketitle
\thispagestyle{empty}



\section{A Rabbit Hunting Problem}  Imagine an invisible rabbit that starts at some unknown integer point \( A \) on the number line. At each time step, it hops by a fixed but unknown integer stride \( B \). Both \( A \) and \( B \) are fixed integers, but their values are unknown. Suppose you have a magic hammer that you can throw at any integer point on the number line at each time step. When the hammer strikes the rabbit, it instantly squeals, indicating you have hit it. The problem now is to devise a strategy that guarantees your  hammer will hit the rabbit in finitely many steps. 

 The above problem is similar to search problems in theoretical computer science, where the goal is to systematically locate a hidden target (in this case, the rabbit). We will provide two algorithms to solve this problem. The first involves Cantor's diagonal trick from set theory, and the second is a probabilistic approach. After presenting both algorithms, we will discuss generalizations showing how our two methods differ. Finally, we end by posing further questions for the reader to investigate.

We will use basic results from set theory, calculus, and probability theory found in standard undergraduate textbooks on these topics; see \cite{apostol1991calculus, apostol2000calculus} for instance. 

\section{A set-theoretic algorithm}

\noindent 
The sequence of hops of this invisible rabbit is determined uniquely by an ordered pair $(A, B)$ of integers. As $\mathbb{Z}^2$  has a countable infinite cardinality, we can enumerate all possible sequences of hops as follows. Let  $r_{ij}$ denote the rabbit's position at  $j$th time step when it hops along the $i$th sequence. 
\begin{eqnarray*}
    r_{11}, r_{12}, r_{13}, \cdots, \\
    r_{21}, r_{22}, r_{23}, \cdots, \\
    r_{31}, r_{32}, r_{33}, \cdots, \\
    r_{41}, r_{42}, r_{43}, \cdots, \\
    \vdots\\
    r_{k1}, r_{22}, r_{k3}, \cdots, \\
    \vdots
\end{eqnarray*}

The hammering strategy is now clear: just hit along the diagonal!
\[r_{11}, r_{22}, r_{33}, r_{44}, \cdots, r_{kk}, \cdots\]
Sooner or later, the hammer is guaranteed to hit our invisible rabbit. To see this, let us suppose the path of our invisible rabbit is $r_{m1}, r_{m2}, r_{m3}, \cdots, r_{mm}, \cdots$ ($m$ is unknown). Then, our hammering sequence will meet this sequence of hops of the rabbit at the $m$th step: $r_{mm}$. 

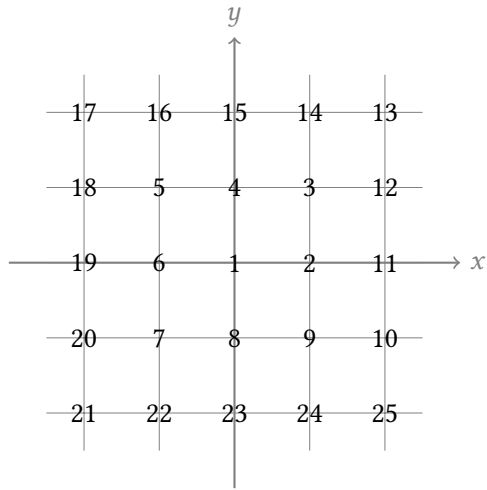
\begin{figure} \label{snake}

\begin{tikzpicture}

\draw[step=1cm, gray, very thin] (-2.5,-2.5) grid (2.5,2.5);

\draw[->, gray, thick] (-3,0) -- (3,0) node[right] {$x$};
\draw[->, gray, thick] (0,-3) -- (0,3) node[above] {$y$};

\node at (0,0) {1};
\node at (1,0) {2};
\node at (1,1) {3};
\node at (0,1) {4};
\node at (-1,1) {5};
\node at (-1,0) {6};
\node at (-1,-1) {7};
\node at (0,-1) {8};
\node at (1,-1) {9};
\node at (2,-1) {10};
\node at (2,0) {11};
\node at (2,1) {12};
\node at (2,2) {13};
\node at (1,2) {14};
\node at (0,2) {15};
\node at (-1,2) {16};
\node at (-2,2) {17};
\node at (-2,1) {18};
\node at (-2,0) {19};
\node at (-2,-1) {20};
\node at (-2,-2) {21};
\node at (-1,-2) {22};
\node at (0,-2) {23};
\node at (1,-2) {24};
\node at (2,-2) {25};

\end{tikzpicture}

\caption{A bijection from $\mathbb{Z}^2 \rar \mathbb{N}$.}
\end{figure}

It is worth noting that the above argument not only proves the existence of a strategy but also gives a constructive algorithm, as we can easily write down an explicit bijection from $\mathbb{Z}^2 \rar \mathbb{N}$. For instance, one bijection comes from following a snake pattern on the lattice:   $1 \mapsto (0, 0), 2 \mapsto (1, 0), 3 \mapsto (1, 1), 4 \mapsto (0, 1)$, etc.; see Figure 1. 

\section{A probabilistic algorithm}

The main idea behind a probabilistic algorithm is to find an appropriate (hammering) function $h \colon \mathbb{N} \rar \mathbb{N}$ such that, at each time step $n$,  hitting a hammer at $X_n \in \mathbb{Z}$, where $X_n$ is a discrete random variable with  uniform distribution
\[X_n \sim \text{Unif}\{ -h(n), -h(n)+1, \dotsc, -1, 0, 1,  \dotsc, h(n)-1, h(n) \},\]
is guaranteed to hit the rabbit in finitely many steps with probability 1.

Note that the rabbit's position at the $n$th time step is given by $R_n = A+Bn$. As $R_n$ is changing at a linear rate, we want our function $h(n)$ to be, at least,  
an increasing function such that for any two integers $A$ and $B$,
\[h(n) > A+Bn \; \; \text{ for } n >\!> 0.\]
The above inequality ensures that our interval $[-h(n), h(n)]$ contains the rabbit for sufficiently large values of $n$, say $n \ge n_0$ for some $n_0$ which may depend on $A$ and $B$.

Finally, our function \( h \) must also satisfy the condition that
\[ P(R_1 \ne X_1, R_2 \ne X_2, \dotsc, R_n \ne X_n, \dotsc) = 0. \]
\noindent
Why is this necessary? This condition states that the probability of never matching the rabbit's position (i.e., the probability that our guessed position \( X_n \) will not coincide with the rabbit's position \( R_n \) for any \( n \)) is zero. In other words, it ensures that our probabilistic method will catch the rabbit in a finite number of steps with probability 1.

To this end, let $h$ be an increasing function that satisfies the above two conditions. Then we have the following for any fixed $A, B$ (which fixes $n_0$):

 \begin{eqnarray*}
   P(R_1 \ne X_1, R_2 \ne X_2, \dotsc, R_n \ne X_n, \dotsc)
   & = & \prod_{k =1}^{\infty} P(R_k \ne X_k)  \\
   & \le & \prod_{k =n_0}^{\infty} P(R_k \ne X_k) \\
   & = & \prod_{k =n_0}^{\infty} \frac{2h(k)}{2h(k)+1.}
 \end{eqnarray*}
Note that 
\begin{eqnarray*}
   \prod_{k =n_0}^{\infty} \frac{2h(k)}{2h(k)+1} \rar 0 & \iff &  \sum_{k = n_0}^{\infty} \log \left(\frac{2h(k)}{2h(k)+1} \right) \rar -\infty \\ 
   & \iff & \sum_{k = n_0}^{\infty} \log \left(\frac{2h(k)+1}{2h(k)} \right) \rar \infty \\ 
   & \iff & \sum_{k = n_0}^{\infty} \log \left(1 + \frac{1}{2h(k)} \right) \rar \infty. \\ 
\end{eqnarray*}

 A well-known fact in the theory of infinite series states that for a sequence $(a_k)$ of positive terms, the series $\sum \log(1 +a_k)$ diverges if and only if $\sum a_k$ diverges.
This gives us the final condition we must impose on our function: $\sum 1/h(k)$ must diverge. 

To summarize, our probabilistic algorithm is guaranteed to hit the rabbit in finitely many steps provided there exists a function $h \colon \mathbb{N} \rar \mathbb{N}$ such that
\begin{enumerate}
    \item[(1)] $h$ is an increasing function,
    \item[(2)] For any two integers $A$ and $B$, $h(k) > A+Bk \; \; \text{ for } k >\!> 0$, and 
    \item[(3)] $\sum \frac{1}{h(k)}$ diverges.
\end{enumerate}

\noindent

To meet these requirements, we need a function \( h(k) \) that increases faster than linearly but slower than quadratically. Specifically, the function's growth rate should be slower than quadratic because the series \(\sum \frac{1}{k^2}\) converges. A natural candidate to try is \( h(k) = k^{1 + \epsilon} \) for \( 0 < \epsilon < 1 \). However, this does not satisfy the requirement (3) because the series \(\sum \frac{1}{k^{1 + \epsilon}}\) converges, as confirmed by the \(p\)-series test.
This means we need a function that grows more slowly than \( k^{1+\epsilon} = kk^{\epsilon}\) for any \( \epsilon > 0 \). A suitable function with this property is \(k\log k\). This suggests examining \( h(k) = \lfloor k \log (k) \rfloor \).

It is clear that this function meets the first two conditions. To verify whether the series \(\sum_{k=2}^{\infty} \frac{1}{\lfloor k \log (k) \rfloor}\) diverges, we can use the integral test. Consider the function \( f(x) = \frac{1}{x \log(x)} \), which is continuous, positive, and decreasing for \( x > 1 \). Evaluating the integral:

\[ 
\int_1^\infty \frac{dx}{x \log x} = \log( \log(x)) \Big|_{1}^{\infty} \rightarrow \infty,
\]
the integral test confirms that \(\sum_{k=2}^{\infty} \frac{1}{k \log(k)}\) diverges. Consequently, using the comparison test, we conclude that \(\sum_{k=2}^{\infty} \frac{1}{\lfloor k \log (k) \rfloor}\) also diverges.
Thus, \( h(k) = \lfloor k \log (k) \rfloor \) satisfies all the required properties for our algorithm.

This gives the desired probabilistic algorithm to find the rabbit in finitely many steps with probability 1  : at time step $n$, hammer at $X_n$ where 
$X_n$ follows a uniform distribution
\[X_n \sim \text{Unif}\{ -\lfloor n \log (n) \rfloor , -\lfloor n \log (n) \rfloor +1, \dotsc, -1, 0, 1,  \dotsc, -\lfloor n \log (n) \rfloor -1, -\lfloor n \log (n) \rfloor  \}.\]
It is worth noting that many functions satisfy conditions (1-3) that are not just a multiple of $\lfloor  x\log(x)\rfloor $. For instance,   $ h(x) =\lfloor x \log(\log(x+1)) \rfloor $ is yet another candidate. 

\subsection{Expected number of steps}

Even though our probabilistic method finds the rabbit in finitely many steps with probability 1, we will now show that the expected number of steps is infinite.

To demonstrate this, let $T$ be the random variable representing the number of steps needed to hit the rabbit using our probabilistic approach. Formally, 
\[ T := \text{min}\{ n \colon X_n = R_n \}. \]
Note that $E[T] = \sum_{k=0}^{\infty} P(T> k)$.  We set $a_k =  P(T> k)$ for simplicity. It is clear from our algorithm that $a_0= a_1 = \cdots a_{n_0-1}=1$. And, for $k \ge n_0$, we have the following.
\begin{eqnarray*}
a_k  &:=& P(T>k) \\
& = & P(\text{rabbit not hit up to time } k) \\
& = & \prod_{n=1}^k  P(\text{rabbit  not hit at time } n) \\
& = & \prod_{n=n_0}^k \left( 1 - \frac{1}{2h(n)}\right) \\
& = & a_{k-1} \left( 1 - \frac{1}{2\lfloor k \log k \rfloor} \right).
\end{eqnarray*}
Since $E(T) = \sum a_k$, we will be done if we can prove that series $\sum a_k$  diverges. To this end, it is natural to apply the ratio test. Unfortunately,
\[ 
\lim_{k \rightarrow \infty} \frac{a_k}{a_{k-1}} 
=  \lim_{k \rightarrow \infty} \left( 1 - \frac{1}{2\lfloor k \log k \rfloor} \right) 
= 1, 
\]
so the ratio test is inconclusive. Instead, we use the following version of Raabe's Test.
\vskip 3mm
\noindent 
\textbf{Raabe's Test:} Let $\sum_{n=1}^{\infty} c_n$ be a series of positive terms. Set $\rho_n := n(c_n/c_{n+1} - 1)$ for all $n \ge 1$.
\begin{itemize} 
\item If $\lim_{n \rightarrow \infty}\rho_n < 1$, then $\sum_{n=1}^{\infty} c_n$ diverges.
\item If $\lim_{n \rightarrow \infty} \rho_n > 1$, then  $\sum_{n=1}^{\infty} c_n$ converges.
\item The test is inconclusive when the limit is equal to 1.
\end{itemize}
\noindent
For our problem, we have the following for all $n \ge n_0$.
\begin{eqnarray*}
\rho_n &=&   n\left ( \frac{a_n}{a_{n+1}} - 1 \right) \\
& = &    n\left ( \left(1 - \frac{1}{2\lfloor (n+1) \log(n+1)\rfloor }\right)^{-1}- 1 \right) \\
& = & n \left ( \frac{2 \lfloor (n+1)\log(n+1) \rfloor}{2\lfloor (n+1) \log (n+1) \rfloor - 1}- 1 \right) \\
& = & \frac{n}{2 \lfloor (n+1) \log (n+1) \rfloor - 1} \rightarrow 0 \; \text{ as } n \rightarrow \infty.
\end{eqnarray*}
Since the limit is less than 1,  we can apply Raabe's test to conclude that the series $\sum a_n$ diverges. In other words, $E(T)= \infty$.

\section{Generalizations}

It is interesting to note that although both approaches solve the problem as stated, they are very different in nature. This is clear if we consider two simple generalizations of the problem.

\begin{enumerate}

\item Suppose the position of the rabbit at time $n$ is given by $R_n = A + Bn + Cn^2$, where $A, B, C \in \mathbb{Z}$ are arbitrary unknown integers. The set-theoretic approach will still work, as the number of possible sequences (which has the same cardinality as $\mathbb{Z}^3$) is still countable. However, the probabilistic approach will no longer work, as $h(n)$ must now grow at least at a quadratic rate to ensure that it eventually overtakes the rabbit for all possible sequences, and this, in turn, implies that $\sum \frac{1}{h(k)}$ does not diverge.

\item Suppose the position of the rabbit at time $n$ is still given by $R_n = A + Bn$, but now $A, B \in \mathbb{R}$ are arbitrary real numbers. The hammer, still thrown at integer points, is considered to hit the rabbit at time $n$ if $\lvert R_n - h(n) \rvert \leq \frac12$. The set-theoretic approach does not work as there are now uncountably many possible sequences along which the rabbit can hop. The probabilistic approach, however, still works (a formal proof is left as an exercise).

\end{enumerate}

\section{Further Questions}

We end with a few questions for the reader to explore.
\begin{enumerate}
\item  The expected number of steps in our probabilistic algorithm was infinite.  Is it possible to devise a different probabilistic algorithm with a finite expected number?
    \item  The rabbit hops we discussed were one-dimensional. One can consider higher-dimensional analogs of this problem. For instance, in two dimensions, imagine a rabbit that starts at some unknown lattice point $(a_1, a_2)$, and at each time step, it hops by a fixed unknown vector $(b_1, b_2)$. The position of the rabbit at the $n$th time step is then given by $R_n = (a_1 + nb_1, a_2 +nb_2)$.  Is it possible to devise a strategy to find the rabbit in finitely many steps?
\item   What properties of the group $(\mathbb{Z}, +)$ played a role in our analysis? Identify those properties and generalize these algorithms to rabbit hops on abstract groups. 
\end{enumerate}

\begingroup
\raggedright

\bibliographystyle{alpha}

\endgroup

\end{document}